\documentclass[11pt,leqno]{article}

\usepackage{amssymb}
\usepackage{amsthm}
\usepackage{eucal}
\usepackage{amsmath,amsfonts}
\usepackage{amsxtra}

\newcommand\Cinf{\mathcal{C}^\infty}

\newcommand{\const}{\mathop{\rm const\,}\nolimits}
\newcommand{\supp}{\mathop{\rm supp\,}\nolimits}

\newcommand\cD{\mathcal D}

\newcommand{\HH}{\mathbf{H}}
\newcommand{\RR}{\mathbf{R}}

\def \tilde{\widetilde}

\def \hat{\widehat}
\newcommand{\st}[1]{\ensuremath{^{\scriptstyle \textrm{#1}}}}

\makeatletter
\@addtoreset{equation}{section}
\makeatother

\newtheorem{thm}{Theorem}[section]
\newtheorem{theorem}[thm]{Theorem}
\newtheorem*{theorem*}{Theorem}
\newtheorem{corollary}[thm]{Corollary}

\newtheorem{lemma}[thm]{Lemma}
\newtheorem*{lemma*}{Lemma}
\newtheorem{proposition}[thm]{Proposition}

\numberwithin{equation}{section}

\newcommand{\romanlist}{
  \renewcommand{\theenumi}{\roman{enumi}}%
  \renewcommand{\labelenumi}{(\theenumi)}%
}

\renewenvironment{proof}[1][Proof]
           {\medbreak\noindent \emph{#1: \enspace}}

\newenvironment{remark}[1][Remark.]
           {\medbreak\noindent \textbf{#1 \enspace }}
           {\par \medbreak}

\makeatletter
\def\iint{\DOTSI\protect\ints@\tw@}
\def\iiint{\DOTSI\protect\ints@\thr@@}
\def\iiiint{\DOTSI\protect\ints@{4}}
\def\idotsint{\DOTSI\protect\ints@\z@}

\def\intkern@{\mkern-6mu\mathchoice{\mkern-3mu}{}{}{}}
\let\DOTSI\relax
\let\ilimits@\displaylimits

\def\ints@#1{%
  \mkern-7mu\mathchoice{\mkern-2mu}{}{}{}%
  \mathop{\mkern7mu\mathchoice{\mkern2mu}{}{}{}%
    \intop\ifnum#1=\z@\intdots@
    \else\intkern@\fi
    \ifnum#1>\tw@\intop\intkern@\fi
    \ifnum#1>\thr@@\intop\intkern@\fi
    \intop
  }\ilimits@
}
\makeatother

\makeatletter
\def\@maketitle{\newpage
 \null
 \vskip 2em
 \begin{center}%
  {\Large\bf \@title \par}%
  \vskip 1.5em
  {\normalsize
   \lineskip .5em
   \begin{tabular}[t]{c}\@author
   \end{tabular}\par}%
  \vskip 2em
  {\@date}%
 \end{center}%
 \par
 \vskip 2.5em}
\makeatother

\begin{document}

\begin{center}
\textbf{Support Theorems for Horocycles on Hyperbolic
  Spaces.}\\[2ex]
Sigurdur Helgason\\[2ex]
In  Memory of S.S. Chern

\end{center}

\vspace{5ex}

\section{Introduction}
\label{sec:1}

The Radon transform associates to a function on a  space $X$ a
function $\hat{f}$ on a family $\Xi$ of subsets $\xi \subset X$
with the definition,
\begin{equation}
  \label{eq:1.1}
  \hat{f} (\xi) = \int_\xi f (x) \, dm (x) \, , \quad \xi \in \Xi\, ,
\end{equation}
$dm$ being a given measure on each $\xi$.  Radon's original
question [9] was whether this mapping  $f \to \hat{f}$ was
injective, in other words whether $f$ is determined by the
integrals (\ref{eq:1.1}).  Along with this injectivity problem,
determining the range of the mapping $f \to \hat{f}$ is an
interesting question.  A part of this question is the so-called
{\em support theorem}.  While the implication

\begin{equation}
  \label{eq:1.2}
  \supp (f) \textup{\,\, compact\,\,} \Rightarrow \supp (\hat{f})
   \textrm{\,\, compact\,\,}
\end{equation}
($\supp$ denoting support) will usually hold for simple reasons,
the converse implication

\begin{equation}
  \label{eq:1.3}
    \supp (\hat{f})\,\, \textrm{compact}\,\, \Rightarrow  \supp (f)
   \textrm{\,\, compact\,\,}
\end{equation}
is designated the support theorem (usually with extra assumption
on $f$).  Positive answers for some examples lead to various
applications:

\romanlist
\begin{enumerate}
\item An explicit description of the range $\cD (X)\sphat$
 where $X$ is a Euclidean space or a symmetric space
  of the noncompact type ([2], [3]). Here $ (\cD = \Cinf_c)$. In the first case,
  $\hat{f}$ in (\ref{eq:1.1}) is integration over hyperplanes  in
  $X =\RR^n$; in the latter case $\hat{f}$ in (\ref{eq:1.1})
  refers to integration over horocycles $\xi$ in the symmetric
  space $X$.

\item  Medical application in X-ray reconstruction ([6], p.47).

\item  Existence theorem for invariant differential equations on
  a symmetric space $X$ ([3], Lemma 8.1 and Theorem 8.2).

\end{enumerate}

While these results rely on special methods for each case,
microlocal analysis has been used e.g.~by Quinto [8] for results
of more general nature, requiring however stronger  {\em a~priori} assumptions about $f$ and its support.

For a symmetric space $X$ of the noncompact type there are two
natural Radon transforms, the X-ray transform and the horocycle
transform; in both cases (\ref{eq:1.3}) holds ([4], [3]).  If
$X$ has rank one, then a horocycle has codimension one and its
interior is well defined ([1]).  Thus one can raise the question
of a support theorem for the X-ray transform $f \to \hat{f}$
relative to a fixed horocycle.  If $f$ is assumed exponentially
decreasing, the support theorem does indeed hold ([5]).
Specifically, a function on $X$ is said to be {\em exponentially
  decreasing} if 
\begin{displaymath}
  \sup_x f (x) e^{m\, d(0,x)}< \infty
\end{displaymath}
for each $m>0, 0 \in X$ denoting the origin and $d$ the distance.

For  $X$ a hyperbolic space we consider in this note the
analogous question for the horocycle transform $f \to \hat{f}$,
relative to a fixed horocycle (Theorem~2.2), extending a
result by Lax and Phillips ([7]).

This paper is dedicated to the memory of Shiing-Shen Chern in
appreciation of his generosity and thoughtfulness through many
years.  In 1959 when I was planning my 1962 book, his
encouragement and advice were invaluable.  He took active
interest in my work developing Radon transform theory for
homogeneous spaces and told me how his 1942 incidence definition
for a pair of homogeneous spaces fits into the program.  In fact,
some of Radon's old results from 1917 are best understood from
this point of view.

\section{The horocycle transform on $H^n$}
\label{sec:2}

For the support question we take the hyperbolic space $H^n$ with the metric
\begin{equation}
  \label{eq:2.1}
  ds^2 = \frac{dx^2_1 + \cdots + dx^2_n}{x^2_n} \, , \quad x_n>0\,.
\end{equation}
%

In the metric (2.1) the
geodesics are the circular arcs perpendicular to the plane $x_n
=0$; among these are the half lines perpendicular to $x_n=0$.
The horocycles perpendicular to these last geodesics are the
planes $x_n =\const$.  The other horocycles are the Euclidean
$(n-1)$-spheres tangential to the boundary.

Let $\xi \subset \HH^n$ be a horocycle in the half space model.
It is a Euclidean sphere with center~$(x',r)$ (where $x' =
(x_1,\ldots x_{n-1})$) and radius~$r$.  We consider the
intersection of $\xi$ with the $x_{n-1}x_n$ plane.  It is the
circle $\gamma : x_{n-1} = r \sin \theta$, $x_n = r (1-\cos
\theta)$ where $\theta$ is the angle measured from the point of
contact of~$\xi$ with $x_n=0$.  The plane $x_n =r (1-\cos
\theta)$ intersects $\xi$ in an $(n-2)$-sphere whose points are
$x'+ r\sin \theta \omega'$ where $\omega' = (\omega_1 , \ldots ,
\omega_{n-1})$ is a point on the unit sphere $S_{n-2}$ in
$\RR^{n-1}$.  Let $d \omega'$ be the surface element on
$S_{n-2}$.

\begin{proposition}
  \label{prop:2.2}
Let $f$ be exponentially decreasing on $\HH^n$.  Then in the
notation above,
\begin{equation}
  \label{eq:2.2}
\hat{f} (\xi)   = \int^\pi_0 \int_{S_{n-2}} f (x'+ r\sin\theta\omega'
    \, ,\,  r (1-\cos \theta))\, d\omega' 
      \biggl( \frac{\sin \theta}{1-\cos \theta}\biggr)^{n-2}\, 
        \frac{d\theta }{1-\cos \theta}\, .
\end{equation}
\end{proposition}

\begin{proof}
%
%
Since horizontal translations preserve (\ref{eq:2.1}) and commute
with $f \to \hat{f}$ we may assume $x'=0$.

The plane $\pi_\theta : x_n = r (1-\cos \theta)$ has the
non-Euclidean metric
\begin{displaymath}
  \frac{dx^2_1 + \cdots + dx^2_{n-1}}{r^2 (1-\cos \theta)^2}
\end{displaymath}
and the intersection $\pi_\theta \cap \xi$ is an
$(n-2)$-sphere with induced metric
\begin{displaymath}
  \frac{r^2 \sin^2 \theta (d\omega')^2}{r^2 (1-\cos \theta)^2}\, ,
\end{displaymath}
where $(d\omega')^2$ is the metric on the $(n-2)$-dimensional
unit sphere in $\RR^{n-1}$.  The non-Euclidean volume element on
$\xi \cap \pi_\theta$ is thus
\begin{displaymath}
  \biggl( \frac{\sin \theta}{1-\cos \theta}\biggr)^{n-2}\,
  d\omega'\, .
\end{displaymath}
The non-Euclidean arc element on $\gamma$ is by (2.1) equal to
$d\theta / (1-\cos \theta)$.  Putting these facts together (2.2)
follows by integrating over $\xi$ by slices $\xi \cap \pi_\theta$.

\end{proof}

\begin{theorem}
  \label{2.2}

Let $\xi_0 \subset \HH^n$ be a fixed horocycle.  Let $f$ be
exponentially decreasing and assume
\begin{displaymath}
  \hat{f} (\xi) =0
\end{displaymath}
for each horocycle $\xi$ lying outside $\xi_0$.  Then 
\begin{displaymath}
  f(x) =0 \hbox{\,\, for \,\,} x \hbox{\,\, outside \,\,}\xi_0\,.
\end{displaymath}

\end{theorem}

\begin{remark}
  
For the case $n=3$ this is proved in Lax-Phillips [7].  As we see
below, this case is an exception and the general case requires
additional methods.

\end{remark}



\begin{proof}
By homogeneity we may take $\xi_0$ as the plane
$x_n=1$.  Assuming $\hat{f} (\xi) =0$ we take the Fourier
transform in the $x'$ variable of the right hand side of (2.2),
in other words integrate it against $e^{-i\langle x',\eta'
  \rangle}$ where $\eta' \in \RR^{n-1}$.

\end{proof}

Then 
\begin{displaymath}
  \int^\pi_0 \int_{S_{n-2}} \tilde{f}  (\eta',r (1-\cos\theta))
     e^{-ir\sin\theta \langle \eta',\omega'\rangle} \, d\omega'
     \biggl( \frac{\sin\theta}{1-cos \theta}  \biggr)^{n-2} \,
     \frac{d\theta}{1-\cos \theta} =0 \, .
\end{displaymath}
By rotational invariance the $\omega'$ integral only depends on
the norm $|\eta'| r \sin \theta$ so we write
\begin{displaymath}
  J (r\sin \theta |\eta'|) = \int_{S_{n-2}}
     e^{-ir\sin\theta \langle \eta',\omega' \rangle} \, d\omega'\,.
\end{displaymath}
and thus
\begin{displaymath}
  \int^\pi_0 \tilde{f} (\eta', r (1-\cos \theta)) 
     J (r\sin\theta |\eta' | )
     \biggl( \frac{\sin\theta}{1-\cos\theta}\biggr)^{n-2}
       \, \frac{d\theta}{1-\cos \theta} =0 \, .
\end{displaymath}
Here we substitute $u=r (1-\cos \theta)$ and obtain
\begin{equation}
  \label{eq:2.3}
  \int^{2r}_0 \tilde{f} (\eta',u) J ((2ur-u^2)^{1/2} |\eta'|)
    \frac{r}{u^{n-1}} (2ur-u^2)^{\frac12 (n-3)}\, du=0\,.
\end{equation}
Since the distance from the origin $(0,1)$ to $(x',u)$ satisfies
\begin{displaymath}
  d ((0,1),(x',u)) \geq d ((0,1),(0,u)) = \int^1_u
    \frac{dx_n}{x_n} = - \log u
\end{displaymath}
so
\begin{displaymath}
  e^{d ((0,1),(x',u))} \geq \frac1u \, ,
\end{displaymath}
and since
\begin{displaymath}
  \tilde{f} (\eta' ,u) = \int_{\RR^{n-1}} f (x',u)
    e^{-i \langle x',\eta' \rangle }\, dx' \, ,
\end{displaymath}
we see from the exponential decrease of~$f$, that the function $u \to
\tilde{f} (\eta',u)/u^{n-1}$ is continuous down to~$u=0$.

\vspace{2ex}

\noindent{{\bf The case $n=3$.}}\quad  In this simplest case (2.3)
takes the form
\begin{equation}
  \label{eq:2.4}
  \int^{2r}_0 \tilde{f} (\eta',u) u^{-2} J ((2ur-u^2)^{1/2}
     |\eta'|) \, du =0 \, .
\end{equation}
We need here standard result for Volterra\index{Volterra integral
equation} integral equation
(cf.~Yosida [10]).

\begin{proposition}
  \label{prop:2.3}

Let $a<b$ and $f \in C [a,b]$ and $K(s,t)$ of class $C^1$ on
$[a,b] \times [a,b]$.  Then the integral equation
\begin{equation}
  \label{eq:2.5}
  \varphi (s) + \int^s_a K (s,t) \varphi (t) \, dt = f (s)
\end{equation}
has a unique continuous solution $\varphi (t)$.  In particular,
if $f \equiv 0$ then $\varphi \equiv 0$.
\end{proposition}

\begin{corollary}
  \label{cor:2.4}
Assume $K (s,s)\neq 0$ for $s \in [a,b]$.  Then the equation
\begin{equation}
  \label{eq:2.6}
  \int^s_a K(s,t) \psi (t) \, dt =0 \hbox{\,\, implies\,\,}
     \psi \equiv 0 \, .
\end{equation}

\end{corollary}

This follows from Prop.~2.3 by differentiation.  Using Cor.~2.4
on (2.4) we deduce $\tilde{f} (\eta',u)=0$ for $u \leq 2r$ with
$2r \leq 1$ proving Theorem~2.2 for $n=3$.

\vspace{2ex}

\noindent{{\bf The case $n=2$.}}\quad Here (2.3) leads to the
generalized Abel integral equation $(0< \alpha <1)$.
\begin{equation}
  \label{eq:2.7}
  \int^s_a \frac{G (s,t)}{(s-t)^\alpha} \varphi (t) \, dt = f (s)\,.
\end{equation}

\begin{theorem}
  \label{th:2.5}

With $f$ continuous, $G$ of class $C^1$ and $G (s,s) \neq 0$ for
all $s \in [a,b]$, equation (2.7) has a unique continuous
solution~$\varphi$.  In particular, $f \equiv 0 \Rightarrow
\varphi \equiv 0$.

\end{theorem}

This is proved by integrating the equation  against
$1/(x-s)^{1-\alpha}$ whereby the statement is reduced to
Cor.~2.4 (cf. Yosida, {\it loc.cit.}).

This proves Theorem 2.2 for $n=2$.

\vspace{2ex}

\noindent{{\bf The general case.}}\quad  Here the parity of $n$
makes a difference.  For $n$ odd we just use the following lemma.

\begin{lemma}
  \label{lem:2.6}
Assume $\varphi = C^1 ([a,b])$ and that $K (s,t)$ has all
derivatives with respect to $s$ up to order $m-2$ equal to $0$
on the diagonal $(s,s)$.  Assume the $(m-1)$\st{th} order
derivative is nowhere $0$ on the diagonal.  Then (2.6) still
holds.

\end{lemma}

In fact, by repeated differentiation of (2.6) one can show that
(2.5) holds with a kernel
\begin{displaymath}
  \frac{K^{(m)} (s,t)}{\{ K^{(m-1)}_s (s,t)\}_{t=s}}
\end{displaymath}
and $f \equiv 0$.

This lemma proves Theorem~2.2 for $n$~odd.  For $n$~even we write
(2.3) in the general form
\begin{equation}
  \label{eq:2.8}
  \int^s_0 F(u) 
H ((su-u^2)^{1/2})(su-u^2)^{\frac12 (n-3)}
\, du =0\quad
n \hbox{\,\, even\,\,} \geq 2\, ,
\end{equation}
where $H (0) \neq 0$.

\begin{theorem}
  \label{th:2.7}

Assume $F \in C ([0,1])$ satisfies (2.8) for $0 \leq s \leq 1$
and $H \in \Cinf$ arbitrary with $H(0) \neq 0$.  Then $F
\equiv 0$ on $[0,1]$.

\end{theorem}

\begin{proof}

We proceed by induction on $n$, the case $n=2$ being covered by
Theorem~2.5.  We assume the theorem holds for~$n$ and any
function $H$ 
satisfying $H (0) \neq 0$.  We
consider (2.8) with $n$~replaced by $n+2$ and take $d/ds$.  The
result is with $H_1 (x) = H' (x) x + (n-1) H (x)$,
\begin{displaymath}
  \int^s_0 F(u)u H_1 ((su-u^2)^{\frac12}) 
     (su-u^2)^{\frac12 (n-3)}   \, du =0\,.
\end{displaymath}
Since 
 $H_1 (0) \neq 0$ we conclude $F \equiv 0$ by
induction.  This finishes the proof of Theorem~2.2.

\end{proof}

\end{document}